\documentclass[12pt]{amsart}
\usepackage{amscd}
\usepackage[mathscr]{eucal}
\usepackage{amssymb}
\usepackage{latexsym}
\usepackage{amsthm}
\usepackage{amsmath}
\usepackage{graphicx}

\theoremstyle{plain}

\newtheorem{thm}{Theorem}[section]
\newtheorem{cor}{Corollary}[section]
\newtheorem{lem}{Lemma}[section]

\theoremstyle{definition}

\newtheorem{prf}{Proof}[section]

\title{Invariant linear functionals on $L^{\infty}(\mathbb{R}_+)$}
\author{Ryoichi Kunisada}

\address{Faculty of Education and Integrated Arts and Science, Waseda University, Shinjuku-ku, Tokyo 169-8050, Japan}

\email{rkunisada@aoni.waseda.jp}

\date{}

\begin{document}
\maketitle

\begin{abstract}
We consider a continuous version of the classical notion of Banach limits, namely, positive linear functionals on $L^{\infty}(\mathbb{R}_+)$ invariant under translations $f(x) \mapsto f(x+s)$ of $L^{\infty}(\mathbb{R}_+)$ for every $s \ge 0$. We give its characterization in terms of the invariance under the operation of a certain linear mapping on $L^{\infty}(\mathbb{R}_+)$. Applications to summability methods are provided in the last section.
\end{abstract}

\bigskip

\section{Introduction}
For simplicity, we use the term `mean' in place of `normalized positive linear functional' throughout the paper. An underlying Banach space $X$ is always a certain function space and the order in $X$ is such that for any $f \in X$, $f \ge 0$ if and only if $f(x) \ge 0$ everywhere or almost everywhere.  

Let $\mathbb{R}_+$ be the positive half $[0, \infty)$ of the real line $\mathbb{R}$ and $L^{\infty}(\mathbb{R}_+)$ be the Banach space of all real-valued essentially bounded measurable functions on $\mathbb{R}_+$. Let $\mathbb{N}_0$ be the set of non-negative integers and $l_{\infty}$ be the Banach space of all real-valued bounded functions on $\mathbb{N}_0$.
The primary objective of this paper is translation invariant means on $L^{\infty}(\mathbb{R}_+)$. For each $s \ge 0$, we consider the following linear operator:
\[T_s :  L^{\infty}(\mathbb{R}_+) \longrightarrow L^{\infty}(\mathbb{R}_+), \quad (T_sf)(x) = f(x+s). \]
And let $T_s^*$ be its adjoint operator. Then we say that $\varphi \in L^{\infty}(\mathbb{R}_+)^*$, the dual space of $L^{\infty}(\mathbb{R}_+)$, is a T-invariant mean if the following conditions hold:
\begin{enumerate}
\item $f \ge 0$ implies $\varphi(f) \ge 0$,
\item $\varphi(1) = 1$.
\item $T_s^*\varphi = \varphi$ for every $s \ge 0$
\end{enumerate}

Let us $C_{ub}(\mathbb{R}_+)$ be the Banach space of all real-valued uniformly continuous bounded functions on $\mathbb{R}_+$. We also consider T-invariant means on $C_{ub}(\mathbb{R}_+)$, which is easier to handle than those on $L^{\infty}(\mathbb{R}_+)$, defined by simply replacing the word `$L^{\infty}(\mathbb{R}_+)$' by `$C_{ub}(\mathbb{R}_+)$' in the definition of T-invariant means on $L^{\infty}(\mathbb{R}_+)$. In other words, T-invariant menas on $C_{ub}(\mathbb{R}_+)$ is the restrictions of T-invariant means on $L^{\infty}(\mathbb{R}_+)$ to its closed subspace $C_{ub}(\mathbb{R}_+)$. Let us  denote by $\mathcal{M}$ the set of all such means. In fact, T-invariant menas on $C_{ub}(\mathbb{R}_+)$ can be viewed as a continuous counterpart of Banach limits on $l_{\infty}$. Recall that $\varphi \in l_{\infty}^*$ is called a Banach limit if the following conditions hold:
\begin{enumerate}
\item $f \ge 0$ implies $\varphi(f) \ge 0$,
\item $\varphi(1) = 1$.
\item $T^*\varphi = \varphi$,
\end{enumerate}
where $T : l_{\infty} \rightarrow l_{\infty}$ is define by $(Tf)(n) = f(n+1)$ and $T^*$ denotes its adjoint operator. Let us denote the set of all Banach limits by $\mathcal{B}$. Banach limits have been studied by several authors, see for example [1], [2], [4], [7]. An important fact is that each Banach limit can be identified with an invariant measure on a certain discrete flow and in a similar way, as we will see in the following section, each T-invariant mean on $C_{ub}(\mathbb{R}_+)$ can be identified with an invariant measure on a certain continuous flow which is the suspension of the discrete flow. 

Recall that for any mean $\varphi$ on $l_{\infty}$, $\varphi$ is a Banach limit if and only if
\[\varphi(f) \le \lim_{n \to \infty} \limsup_{m \to \infty} \frac{1}{n} \sum_{i=0}^{n-1} f(m+i) \]
holds for every $f \in l_{\infty}$. A similar characterization of T-invarian means on $C_{ub}(\mathbb{R}_+)$ holds; namely, for any mean $\varphi$ on $C_{ub}(\mathbb{R}_+)$, $\varphi$ is a T-invariant mean if and only if
\[\varphi(f) \le \lim_{\theta \to \infty} \limsup_{x \to \infty} \frac{1}{\theta} \int_{x}^{x + \theta} f(t)dt. \]
holds for every $f \in C_{ub}(\mathbb{R}_+)$. Now this result leads us to define a class $\mathcal{M}_1$ of menas $\varphi$ on $L^{\infty}(\mathbb{R}_+)$ which satisfy the above inequality for every $f \in L^{\infty}(\mathbb{R}_+)$. This class $\mathcal{M}_1$ is our main interest of this paper. It is easy to show that each $\varphi \in \mathcal{M}_1$ is T-invariant, thougt in contrast to the case of $C_{ub}(\mathbb{R}_+)$ a T-invariant mean on $L^{\infty}(\mathbb{R}_+)$ need not satisfy this condition. On the other hand, we characterize the class $\mathcal{M}_1$ by the invariance with respect to a certain linear transformation on $L^{\infty}(\mathbb{R}_+)$. We also consider invariant means with respect to the action of the multiplicative group $\mathbb{R}^{\times}$ of $\mathbb{R}$ in place of the additive group $\mathbb{R}$. 

The papaer is organized as follows. Section 2 deal with elementary results concerning T-invariant means on $C_{ub}(\mathbb{R}_+)$, including example of T-invariant means on $C_{ub}(\mathbb{R}_+)$ of a simple form which generate whole $\mathcal{M}$ as its closed convex hull. In Section 3 we deal with the class $\mathcal{M}_1$ of T-invariant means on $L^{\infty}(\mathbb{R}_+)$. Section 4 is deveoted to the study of invariant means under the action of multiplicative group. Definitions and results in this section is similar to those in Section 3. Section 5 contains applications to summability methods.

\medskip
\section{Preliminary results}
Throughout the paper, we will use the notion of the limit along an ultrafilter $\mathcal{U}$, denoted by $\mathcal{U}\mathchar`-\lim$, which is a generalization of the ordinary definitions of limit along a sequence $\lim_{n \to \infty}$ or a continuous parameter $\lim_{x \to \infty}$. We give its definition in the general setting. Let $f : X \rightarrow Y$ be a mapping of a set $X$ into a compact space $Y$ and $\mathcal{U}$ be an ultrafilter on $X$. Then there exists an element $y$ of $Y$ such that $f^{-1}(U) \in \mathcal{U}$ holds for every neighborhood $U$ of $y$. This element $y$ of $Y$ is called the limit of $f$ along $\mathcal{U}$ and denoted by $\mathcal{U}\mathchar`-\lim_x f(x)$.

Since the classical notion of Banach limits has a close relation to T-invariant means on $C_{ub}(\mathbb{R}_+)$, we first give an overview of this notion.  As is well known, Banach limits can be viewed as invariant measures on a discrete flow defined as follows: Let $\beta\mathbb{N}_0$ be the Stone-\v{C}ech compactification of $\mathbb{N}_0$ and let $\mathbb{N}_0^*$ be the growth $\beta\mathbb{N}_0 \setminus \mathbb{N}_0$ of $\beta\mathbb{N}_0$. 
We denote the translation of $\mathbb{N}_0$ by $\tau_0$;
\[\tau_0 : \mathbb{N}_0 \longrightarrow \mathbb{N}_0, \quad \tau_0(n) = n+1. \]
Then we extend it continuously to $\beta\mathbb{N}_0$ and denote it by $\tau$. Restricting $\tau$ to $\mathbb{N}_0^*$, we get a homeomorphism of $\mathbb{N}_0^*$ onto itself;
\[\tau : \mathbb{N}_0^* \longrightarrow \mathbb{N}_0^*. \]
Then the pair $(\mathbb{N}_0^*, \tau)$ is a discrete flow. We denote by $\mathcal{M}_{\tau}^d$ the set of all $\tau$-invariant Borel probability measures on $\mathbb{N}_0^*$. Then it is known that $\mathcal{B} \cong \mathcal{M}_{\tau}^d$ holds.

Now we take up T-invariant means $\mathcal{M}_{\tau}$ on $C_{ub}(\mathbb{R}_+)$. Similarly, one can interpret them as invariant Borel measures on a certain continuous flow. Since $C_{ub}(\mathbb{R}_+)$ is a Banach algebra, there exists a compact space $\Omega$, which is in fact the maximal ideal space of $C_{ub}(\mathbb{R}_+)$, such that $C_{ub}(\mathbb{R}_+)$ is isomorphic to $C(\Omega)$ of the space of all real-valued continuous functions on $\Omega$. 
The construction of $\Omega$ is as folllows (see [5] for details): consdier a product space $\mathbb{N}_0 \times [0, 1]$ and define an equivalent relation $\sim$ on it by $(\tau\eta, 0) \sim (\eta, 1)$ for all $\eta \in \beta\mathbb{N}_0$. Then $\Omega$ is homeomorphic to the quotient space $(\mathbb{N}_0 \times [0, 1]) / \sim$. Since the subspace $\{(n, \eta) : n \in \mathbb{N}_0, t \in [0,1]\}$ of $\Omega$ is homeomorphic to $\mathbb{R}_+$, it is noted that $\Omega$ is a compactification of $\mathbb{R}_+$ to which every uniformly continuos bounded function on $\mathbb{R}_+$ can be extended continuously. 
We denote by $\overline{f} \in C(\Omega)$ the continuous extension of $f \in C_{ub}(\mathbb{R}_+)$ to $\Omega$. Identifying $\omega = (\eta, t) \in \Omega$ with an ultrafilter $\{A+t : A \in \eta \}$ on $\mathbb{R}_+$(Recall that each element of $\beta\mathbb{N}_0$ can be indentified with an ultrafilter on $\mathbb{N}_0$), $\overline{f}(\omega)$ is given by the formula
\[\overline{f}(\omega) = \omega\mathchar`-\lim_x f(x). \]

Therefore, every mean $\varphi$ on $C_{ub}(\mathbb{R}_+)$ can be identified with a mean on $C(\Omega)$. Thus, by the Riesz representation theorem, there exists a Borel probability measure $\mu$ on $\Omega$ such that
\[\varphi(f) = \int_{\Omega} \overline{f}(\omega) d\mu(\omega) \]
holds for every $f \in C_{ub}(\mathbb{R}_+)$.

Next we cosider an extension of the following semi-flow to $\Omega$; 
\[\tau_0^s : \mathbb{R}_+ \longrightarrow \mathbb{R}_+, \quad \tau_0^s x = x + s, \quad s \ge 0. \]
For each $s \ge 0$, we define linear operator $\overline{T}_s$ by
\[\overline{T}_s : C(\Omega) \longrightarrow C(\Omega), \quad \overline{T}_s\overline{f} = \overline{T_sf}. \]
Let $\overline{T}_s^*$ be its adjoint operator. Notice that $\Omega$ can be regarded as a subset of the positive part of the unit sphere $S^+_{C(\Omega)^*}$ of $C(\Omega)^*$, the dual space of $C(\Omega)$. Then we can consider the restriction of $\overline{T}_s^*$ to $\Omega$, denoted by $\tau^s$;
\[\tau^s : \Omega \longrightarrow \Omega, \quad s \ge 0. \]
Then by the above formula of $\overline{f}(\omega)$, we have
\[\tau^s\omega = \tau^s(\eta, t) = (\tau^{[t+s]}\eta, t+s-[t+s]),  \]
where $s \in \mathbb{R}$ and $[x]$ denotes the largest integer not exceeding a real number $x$. In particular, the restriction of each $\tau^s$ to $\Omega^* = \Omega \setminus \mathbb{R}_+$ is a homeomorphism and the pair $(\Omega^*, \{\tau^s\}_{s \in \mathbb{R}})$ is a continuous flow. Then we have
\[\overline{T_sf}(\omega) = \omega(T_sf) = (T_s^*\omega)(f) = (\tau^s\omega)(f) = \overline{f}(\tau^s\omega). \]
If $\varphi \in \mathcal{M}_{\tau}$, i.e., $\varphi(T_sf) = \varphi(f)$ holds for every $f \in C_{ub}(\mathbb{R}_+)$ and $s \ge 0$, we have 
\[\int_{\Omega^*} \overline{f}(\omega) d\mu(\omega)   
= \int_{\Omega^*} \overline{T_sf}(\omega) d\mu(\omega)  = \int_{\Omega*} \overline{f}(\tau^s\omega) d\mu(\omega)
= \int_{\Omega^*} \overline{f}(\omega) d(\overline{T}_s^*\mu)(\omega)  \]
for each $f \in C_{ub}(\mathbb{R}_+)$ and $s \in \mathbb{R}$. Hence if a mean $\varphi$ is T-invariant then the corresponding measure $\mu$ is an invariant measure, i.e. $\mu(\tau^sA) = \tau(A)$ holds for every Borel set $A$ of $\Omega^*$ and $s \in \mathbb{R}$.

Notice that, by the definition, the continuous flow $(\Omega^*, \{\tau^s\}_{s \in \mathbb{R}})$ is the suspension of the discrete flow $(\mathbb{N}_0^*, \tau)$. Thus the following result follows immediately.

\begin{thm}
$\mathcal{M}_{\tau}^d$ is affinely homeomorphic to $\mathcal{M}_{\tau}$.
\end{thm}

Now we give below examples of T-invariant means on $C_{ub}(\mathbb{R}_+)$. 
Remark that for given $f \in C_{ub}(\mathbb{R}_+)$ and $\omega \in \Omega^*$, the function $f_{\omega}(s) := \overline{f}(\tau^s\omega)$ of a real variable $s \in \mathbb{R}$, the restriction of a continuous function $\overline{f}$ on $\Omega^*$ to the orbit of $\omega$, is also a uniformly continuous bounded function on $\mathbb{R}$. Let $\omega \in \Omega^*$ and $\mathcal{U}$ be an ultrafilter on $\mathbb{R}_+$ not containing any bounded set of $\mathbb{R}_+$. Then we define for each $f \in C_{ub}(\mathbb{R}_+)$
\[\varphi^{\mathcal{U}}_{\omega}(f) = \mathcal{U}\mathchar`-\lim_{x} \frac{1}{x} \int_0^x f_{\omega}(t)dt. \]
It is obvious that each $\varphi^{\mathcal{U}}_{\omega}$ is an T-invariant mean. We denote the set of all such T-invariant means by
$\mathcal{Q}$. The following assertion can be regarded as a continuous version of [4, Theorem 3] and the proof is essentially a simplification of the proof of it.

\begin{thm}
$\mathcal{M} = \overline{co}(\mathcal{Q})$, where $\overline{co}(\mathcal{Q})$ represents the closed convex hull of $\mathcal{Q}$.
\end{thm}

\begin{prf}
By the Krein-Milman theorem, it is sufficient to prove that
\[\sup_{\varphi \in ex(\mathcal{M}_{\tau})} \varphi(f) =\sup_{\varphi \in \mathcal{M}_{\tau}} \varphi(f) = \sup_{\varphi^{\mathcal{U}}_{\omega} \in \mathcal{Q}} \varphi^{\mathcal{U}}_{\omega}(f)\]
for every $f \in C_{ub}(\mathbb{R}_+)$, where $ex(\mathcal{M}_{\tau})$ denotes the set of extreme points of $\mathcal{M}_{\tau}$. 
Notice that the corresponding Bore probability measure $\mu$ on $\Omega^*$ of $\varphi \in ex(\mathcal{M}_{\tau})$ is an ergodic measure. Then by Birkhoff's ergodic theorem, for each $f \in C_{ub}(\mathbb{R}_+)$ we have 
\[\lim_{x \to \infty} \frac{1}{x} \int_0^x \overline{f}(\tau^s\omega)ds = \lim_{x \to \infty} \frac{1}{x} \int_0^x f_{\omega}(t)dt = \int_{\Omega^*} f(\omega) d\mu(\omega) = \varphi(f) \]
for all $\omega \in \Omega^*$ except for some points which form a set of $\mu$-measure $0$. For any such a point $\omega$ and any $\mathcal{U}$, $\varphi^{\mathcal{U}}_{\omega}(f) = \varphi(f)$ holds. The assertion follows immediately.
\end{prf}

Next we define a subadditive functional $\overline{M}_1$ on $L^{\infty}(\mathbb{R}_+)$ by
\[\overline{M}_1(f) = \lim_{\theta \to \infty} \limsup_{x \to \infty} \frac{1}{\theta} \int_{x}^{x + \theta} f(t)dt. \]
Then we have the following discription of T-invarian means on $C_{ub}(\mathbb{R}_+)$.

\begin{thm}
For a mean $\varphi$ on $C_{ub}(\mathbb{R}_+)$, $\varphi$ is T-invariant if and only if 
\[\varphi(f) \le \overline{M}_1(f) \]
holds for every $f \in C_{ub}(\mathbb{R}_+)$.
\end{thm}

\begin{prf}
First, we prove the sufficiency.
For any $f \in C_{ub}(\mathbb{R}_+)$ and $s \ge 0$, we have 
\begin{align}
\varphi(f-T_sf) &\le \lim_{\theta \to \infty} \limsup_{x \to \infty} \frac{1}{\theta} \int_{x}^{x + \theta} (f(t)-f(t+s))dt \notag \\   
&= \lim_{\theta \to \infty} \left(\frac{1}{\theta} \int_x^{x+s} f(t)dt - \frac{1}{\theta} \int_{\theta+x}^{\theta+x+s} f(t)dt\right) \notag \\
&\le \lim_{\theta \to \infty} \frac{2s}{\theta}\|f\|_{\infty} = 0 \notag.
\end{align}    
$\varphi(f -T_sf) \ge 0$ can be proved in a similar way. Thus we have $\varphi \in \mathcal{M}$. Next we prove the necessity. 
First, since $f$ is uniformly continous we have 
\[\lim_{s \to 0} \|T_sf - f\|_{\infty} = 0. \]
Then by the continuity and invariance of $\varphi$, we have for each $\theta > 0$,
\[ \varphi(f) = \frac{1}{\theta} \int_0^{\theta} \varphi(T_sf)ds = \varphi\left(\frac{1}{\theta} \int_0^{\theta} T_sfds\right) \le \limsup_{x \to \infty} \frac{1}{\theta} \int_0^{\theta} f(x+s)ds.  \]
Hence we get
\[\varphi(f) \le \lim_{\theta \to \infty} \limsup_{x \to \infty} \frac{1}{\theta} \int_{x}^{x + \theta} f(t)dt = \overline{M}_1(f). \]
The proof is complete.
\end{prf}

\section{Main results for additive group}
In this section we consider T-invariant means on $L^{\infty}(\mathbb{R}_+)$ and extend some of the preceding results to this case. We denote by $\mathcal{M}_1$ the set of means $\varphi$ on $L^{\infty}(\mathbb{R}_+)$ for which
\[\varphi(f) \le \lim_{\theta \to \infty} \limsup_{x \to \infty} \frac{1}{\theta} \int_{x}^{x + \theta} f(t)dt, \]
holds for every $f \in L^{\infty}(\mathbb{R}_+)$. Then in the same way as the proof of sufficiency in Theorem 2.3, it is shown that elements of $\mathcal{M}_1$ are T-ivariant means on $L^{\infty}(\mathbb{R}_+)$.

In the following, we will identify each $\omega = (\eta, t) \in \Omega$ with the ultrafilter $\{t+A : A \in \eta \}$ on $\mathbb{R}_+$. Given $f(x) \in L^{\infty}(\mathbb{R}_+)$, we consider the set of its translates $\{f_s(x)\}_{s \ge 0} \subseteq L^{\infty}(\mathbb{R}_+)$, where $f_s(x) = f(x+s)$. Then notice that this is a bounded set of $L^{\infty}(\mathbb{R}_+)$ and hence is a weak* relatively compact subset of $L^{\infty}(\mathbb{R}_+)$. Thus for any $\omega \in \Omega^*$ we can define its limit 
along $\omega$ with respect to weak* topology of $L^{\infty}(\mathbb{R}_+)$:
\[f_{\omega}(x) = \omega\mathchar`-\lim_s f_s(x). \]
Though the function $f_{\omega}(x)$ is in $L^{\infty}(\mathbb{R}_+)$ by the definition, that is, $f_{\omega}$ is defined on $\mathbb{R}_+$, we can extend it the whole line $\mathbb{R}$ in a natural way as follows.
\[f_{\omega}(x) = f_{\tau^{-N}\omega}(N+x), \quad x \in [-N, 0], \]
for every $N > 0$. In this way, we consider $f_{\omega}$ to be a function defined on $\mathbb{R}$, that is, in $L^{\infty}(\mathbb{R})$.
Remark that for a function $f(x) \in C_{ub}(\mathbb{R}_+)$, $f_{\omega}(x)$ is equal to the one defined in the previous section. An important fact concerning to this notion is the lemma below. Let us define a subalgebra $\mathfrak{U}$ of $C_{ub}(\mathbb{R}_+)$ as
\[\mathfrak{U} = \{f(x) \in C_{ub}(\mathbb{R}_+) : f^{\prime}(x) \in L^{\infty}(\mathbb{R}_+) \}, \]
where $f^{\prime}$ is the derivative of $f$. In other words, $f(x)$ is in $\mathfrak{U}$ if and only if $f(x)$ is a bounded Lipschitz continuous function on $\mathbb{R}_+$; namely, $f(x)$ is a bounded function on $\mathbb{R}_+$ such that
\[|f(x) - f(y)| \le K|x-y| \]
holds for every pair $x, y$ of $\mathbb{R}_+$, where $K>0$ is some constant.
It is easy to see that if $f(x)$ is in $\mathfrak{U}$ then $f_{\omega}(x)$ is also a bounded Lipschitz continuous function on $\mathbb{R}$ for any $\omega \in \Omega^*$ and hence the derivative $(f_{\omega})^{\prime}(x)$ exists and bounded a.e on $\mathbb{R}$. 

\begin{lem} 
Let $f \in \mathfrak{U}$ and $f^{\prime}$ be its derivative. Then $(f_{\omega})^{\prime}(x) = (f^{\prime})_{\omega}(x)$ holds.
\end{lem}

\begin{prf}
By the definition of $f_{\omega}$,
\[(f^{\prime})_{\omega}(x) = \omega\mathchar`-\lim_s f^{\prime}(x+s). \] 
Then for every $x \ge 0$, 
\[\int_0^x (f^{\prime})_{\omega}(t)dt = \omega\mathchar`-\lim_s \int_0^x f^{\prime}(t+s)dt = \omega\mathchar`-\lim_s (f(x+s) - f(s)) = f_{\omega}(x) - f_{\omega}(0).  \]
Therefore, we have
\[(f^{\prime})_{\omega}(x) = (f_{\omega})^{\prime}(x). \]
\end{prf}

\medskip

For any $\theta > 0$, we define the linear operator $U_{\theta}$ by
\[U_{\theta} : L^{\infty}(\mathbb{R}_+) \longrightarrow C_{ub}(\mathbb{R}_+), \quad (U_{\theta}f)(x) = \frac{1}{\theta} \int_x^{x+\theta} f(t)dt. \] 

\begin{lem}
For any $f(x) \in L^{\infty}(\mathbb{R}_+), \omega \in \Omega^*$ and $\theta > 0$, $(U_{\theta}f)_{\omega}(x) = (U_{\theta}f_{\omega})(x)$ holds for every $x \in \mathbb{R}$. 
\end{lem}

\begin{prf}
By the definition, we have
\[(U_{\theta}f)_{\omega}(x) =\tau^x \omega\mathchar`-\lim_s \frac{1}{\theta} \int_s^{s+{\theta}} f(t)dt =\tau^x \omega\mathchar`-\lim_s \frac{1}{\theta} \int_0^{\theta} f_s(t)dt. \]
Then since $f_s \rightarrow f_{\omega}$ is weak* convergence, the right side of the equation is equal to
\[\frac{1}{\theta} \int_0^{\theta} f_{\tau^x \omega}(t)dt = \frac{1}{\theta} \int_x^{x+\theta} f_{\omega}(t)dt = (U_{\theta}f_{\omega})(x). \]
This completes the proof.
\end{prf}

We define another class $\mathcal{R}$ of means on $L^{\infty}(\mathbb{R}_+)$ satisfying the following condition:
\[\varphi(f) \le \limsup_{x \to \infty} \frac{1}{e^x}\int_0^x f(t)e^tdt \]
for every $f \in L^{\infty}(\mathbb{R}_+)$. Let us introduce the linear operator $S$ defined by
\[S : L^{\infty}(\mathbb{R}_+) \longrightarrow L^{\infty}(\mathbb{R}_+), \quad (Sf)(x) = \frac{1}{e^x}\int_0^x f(t)e^tdt.  \]

\begin{lem}
For any $f \in L^{\infty}(\mathbb{R}_+)$, $(Sf)(x)$ is in $\mathfrak{U}$.
\end{lem}
 
\begin{prf}
A direct computation shows that
\[(Sf)^{\prime}(x) = f(x) - (Sf)(x) , \quad x \ge 0, \]
which gives the result since the right side is bounded.
\end{prf}

We give the converse of this result. Let us define
\[L_0^{\infty}(\mathbb{R}_+) = \{f(x) \in L^{\infty}(\mathbb{R}_+) : f(x) \rightarrow 0 \ as \ x \to \infty \}.    \]

\begin{lem}
Every function $f(x)$ in $\mathfrak{U}$ can be written as $(Sf)(x) + h(x)$ for some $f(x) \in L^{\infty}(\mathbb{R}_+)$ and $h(x) \in L^{\infty}_0(\mathbb{R}_+)$.
\end{lem}

\begin{prf}
By the assumption, $\xi(x) = f(x) + f^{\prime}(x)$ is in $L^{\infty}(\mathbb{R}_+)$ and we have 
\begin{align}
e^x \cdot (f(x) + f^{\prime}(x)) = e^x \cdot \xi(x) &\Longleftrightarrow 
(e^x \cdot f(x))^{\prime} = e^x \cdot \xi(x) \notag \\
&\Longleftrightarrow e^x \cdot f(x) - f(0) = \int_0^x \xi(t) \cdot e^t dt \notag \\
&\Longleftrightarrow f(x) = \frac{1}{e^x}\int_0^x \xi(t)e^tdt + \frac{f(0)}{e^x}, \notag
\end{align}
which proves the theorem.
\end{prf}

We define $\Phi = \{Sf : f \in L^{\infty}(\mathbb{R}_+)\}$ and the above two lemmas shows that $\Phi / (\Phi \cap L^{\infty}_0(\mathbb{R}_+)) = \mathfrak{U} / (\mathfrak{U} \cap L^{\infty}_0(\mathbb{R}_+))$. Also we introduce the two spaces $\mathfrak{U}^{\prime} = \{f^{\prime}(x) : f(x) \in \mathfrak{U}\}$ and $\Phi^{\prime} = \{f-Sf : f \in L^{\infty}(\mathbb{R}_+)\}$.
Then by the proofs of the above two lemmas, it also holds that $\Phi^{\prime} / (\Phi^{\prime} \cap L^{\infty}_0(\mathbb{R}_+)) = \mathfrak{U}^{\prime} / (\mathfrak{U}^{\prime} \cap L^{\infty}_0(\mathbb{R}_+))$. 

The following is the simplest examples of elements of $\mathcal{R}$:
\[\chi_{\omega}(f) = \omega\mathchar`-\lim_x \frac{1}{e^x}\int_0^x f(t)e^tdt, \]
where $\omega \in \Omega^*$. We denote by $\tilde{\mathcal{R}}$ the set of all such means. Then we have the following results, which we will prove in the next section.

\begin{thm} 
$ex(\mathcal{R}) = \tilde{\mathcal{R}}$, where $ex(\mathcal{R})$ denotes the set of extreme points of $\mathcal{R}$.
\end{thm}

\begin{thm}
For any $\varphi \in \mathcal{R}$ there exists a unique probability measure $\mu$ on $\Omega^*$ such that
\[\varphi(f) = \int_{\Omega^*} \chi_{\omega}(f) d\mu(\omega) \]
holds for every $f \in L^{\infty}(\mathbb{R}_+)$.
\end{thm}

Next we give another expression of $\chi_{\omega}$ which plays an important role in the remainder of this paper.
\begin{thm}
For every $f \in L^{\infty}({\mathbb{R}_+})$ and $\omega \in \Omega^*$, it holds that
\[\chi_{\omega}(f) = \int_0^{\infty} f_{\omega}(-t)e^{-t}dt. \]
\end{thm}

\begin{prf}
We begin with the equation
\[f(x) = (Sf)(x) + (Sf)^{\prime}(x), \quad x \ge 0. \]
Then by Lemma 3.1 we get for each $\omega \in \Omega^*$
\[f_{\omega}(x) = (Sf)_{\omega}(x) + ((Sf)_{\omega})^{\prime}(x), \quad x \in \mathbb{R}. \]
As is the proof of Lemma 3.4, we have
\[e^x \cdot (Sf)_{\omega}(x) - (Sf)_{\omega}(0) = \int_0^x f_{\omega}(t) \cdot e^t dt, \quad x \in \mathbb{R}. \]
Hence letting $x \to -\infty$, we have 
\[\chi_{\omega}(f) = (Sf)_{\omega}(0) = -\int_0^{-\infty} f_{\omega}(t)e^t dt = \int_0^{\infty} f_{\omega}(-t)e^{-t}dt. \]
\end{prf}

\begin{lem}
For each $\varphi \in \mathcal{R}$ and $\theta > 0$, 
\[\varphi(U_{\theta}f) = \frac{1}{\theta} \int_0^{\theta} \varphi(T_sf)ds. \]
holds for every $f \in L^{\infty}(\mathbb{R}_+)$.
\end{lem}

\begin{prf}
First we prove for the elements of $\tilde{\mathcal{R}}$. By Theorem 3.3 and Lemma 3.2, we have for each $\omega \in \Omega^*$ and $\theta > 0$,
\begin{align}
\chi_{\omega}(U_{\theta}f) &= \int_0^{\infty} (U_{\theta}f)_{\omega}(-t)e^{-t}dt \notag \\
&= \int_0^{\infty} \left(\frac{1}{\theta}\int_{-t}^{-t+\theta} f_{\omega}(s)ds\right)e^{-t}dt \notag \\
&= \int_0^{\infty} \left(\frac{1}{\theta}\int_0^{\theta} f_{\omega}(s-t)ds\right)e^{-t}dt \notag 
\end{align}
\begin{align}
&= \frac{1}{\theta} \int_0^{\theta} \left(\int_0^{\infty} f_{\omega}(s-t)e^{-t}dt\right)ds \notag \\
&= \frac{1}{\theta} \int_0^{\theta} \chi_{\omega}(T_sf)ds. \notag
\end{align}
Next by Theorem 3.2, for each $\varphi \in \mathcal{R}$ there exists some $\mu \in P(\Omega^*)$ such that
\[\varphi(f) = \int_{\Omega^*} \chi_{\omega}(f) d\mu(\omega). \]
Thus we have 
\begin{align}
\varphi(U_{\theta}f) &= \int_{\Omega^*} \chi_{\omega}(U_{\theta}f)d\mu(\omega) = \int_{\Omega^*} \left(\frac{1}{\theta} \int_0^{\theta} \chi_{\omega}(T_sf)ds\right)d\mu(\omega) \notag \\
&= \frac{1}{\theta}\int_0^{\theta}\left(\int_{\Omega^*} \chi_{\omega}(T_sf) d\mu(\omega)\right)ds = \frac{1}{\theta}\int_0^{\theta}\left(\int_{\Omega^*} \chi_{\omega}(T_sf)d\mu(\omega)\right)ds \notag \\
&= \frac{1}{\theta} \int_0^{\theta} \varphi(T_sf)ds. \notag
\end{align}

\end{prf}

\begin{thm}
For any mean $\varphi$ on $L^{\infty}(\mathbb{R}_+)$, $\varphi \in \mathcal{M}_1$ if and only if $\varphi = 0$ on $\Phi^{\prime}$. Namely, $\varphi \in \mathcal{M}_1$ if and only if $\varphi$ is $S$-invariant.
\end{thm}

\begin{prf} 
(Necessity) Notice that it is sufficient to show that if $\varphi \in \mathcal{M}$, then $M_1(f) = 0$ for every $f(x) \in \Phi^{\prime}$. Let us denote $f(x) = g^{\prime}(x)$ for some $g(x) \in \Phi$. Then for any $\theta > 0$ we have
\[\frac{1}{\theta} \int_x^{x+\theta} f(t)dt = \frac{1}{\theta} \int_x^{x+\theta} g^{\prime}(t)dt =\frac{g(x+\theta) - g(x)}{\theta} \]
\[\therefore \left|\frac{1}{\theta} \int_x^{x+\theta} f(t)dt\right| \le \frac{2\|g\|_{\infty}}{\theta}, \]
which shows that $M_1(f) = 0$. This completes the proof.

(Sufficiency) Now suppose that $\varphi = 0$ on $\Phi^{\prime}$. We first show that $\varphi$ is invarian on $C_{ub}(\mathbb{R}_+)$. It means that $\varphi(T_{\theta}f - f) = 0$ for every $f \in C_{ub}({\mathbb{R}}_+)$ and $\theta \ge 0$. For this it is sufficient to show that $\varphi(T_{\theta}f - f) = 0$ for every $f \in \Phi$ and $\theta \ge 0$ since $\Phi$ is a dense subalgebra of $C_{ub}(\mathbb{R}_+)$. In this case, notice that
\[(T_{\theta}f)(x) - f(x) = f(x+\theta) - f(x) = \int_x^{x+\theta} f^{\prime}(t)dt = \theta \cdot (U_{\theta}f^{\prime})(x). \]
Hence it is sufficient to prove that $\varphi(U_{\theta}g) = 0$ for every $g \in \Phi^{\prime}$. Notice that $\varphi \in \mathcal{R}$ by the assumption that $\varphi = 0$ on $\Phi^{\prime}$. In fact, 
\[\varphi(f) = \varphi((Sf) + (Sf)^{\prime}) = \varphi(Sf) \le \limsup_{x \to \infty} (Sf)(x) = \overline{R}(f). \]
Hence by Lemma 3.5 and the observation that $\Phi^{\prime}$ is invariant under $T_x$, it follows that $\varphi(U_{\theta}g) = \frac{1}{\theta}\int_0^{\theta} \varphi(T_xg)dx = 0$. Therefore by Theorem 2.3, we have
\[\varphi(f) \le \overline{M}_1(f) \]
for every $f \in C_{ub}(\mathbb{R}_+)$. For $f \in L^{\infty}(\mathbb{R}_+)$ in general, we have
\[\varphi(f) = \varphi((Sf) + (Sf)^{\prime}) = \varphi(Sf) \le \overline{M}_1(Sf) = \overline{M}_1(f). \]
This shows $\varphi \in \mathcal{M}_1$.
\end{prf}

Remark that we have shown in the above proof the following result.

\begin{thm}
$\overline{M}_1(f) \le \overline{R}(f)$ for every $f \in L^{\infty}(\mathbb{R}_+)$.
\end{thm}

Next coroallary shows that each $\varphi \in \mathcal{M}_1$ is exactly determined by the values on $C_{ub}(\mathbb{R}_+)$.

\begin{cor}
For each $\varphi \in \mathcal{M}_1$ and $\theta > 0$, $\varphi(f) = \varphi(U_{\theta}f)$ holds.
\end{cor}

\begin{prf}
By Lemma 3.5 and Theorem 3.5, we have that if $\varphi \in \mathcal{M}_1$ then $\varphi(U_{\theta}f) = \frac{1}{\theta} \int_0^{\theta} \varphi(T_sf)ds = \frac{1}{\theta} \int_0^{\theta} \varphi(f)ds = \varphi(f)$.
\end{prf}

\begin{thm}
$\mathcal{M}_1$ is affinely homeomorphic to $\mathcal{M}_{\tau}$.
\end{thm}

Hence it is natural to ask that how can one express the extensions of the elements of $\mathcal{Q}$. The answer to this question is given as follows. Let us denote the extension of $\varphi^{\mathcal{U}}_{\omega}$ by $\overline{\varphi}^{\mathcal{U}}_{\omega}$.

\begin{thm}
For every $f \in L^{\infty}(\mathbb{R}_+)$ and $\varphi_{\omega}^{\mathcal{U}} \in \mathcal{Q}$, it holds that
\[\overline{\varphi}^{\mathcal{U}}_{\omega}(f) = \mathcal{U}\mathchar`-\lim_{x} \frac{1}{x} \int_0^x f_{\omega}(t)dt. \]
\end{thm}

\begin{prf}
By Corollary 3.1, notice that $\overline{\varphi}^{\mathcal{U}}_{\omega}(f) = \overline{\varphi}^{\mathcal{U}}_{\omega}(U_{\theta}f) = \varphi^{\mathcal{U}}_{\omega}(U_{\theta}f)$ holds for each $\theta > 0$. Also by Lemma 3.2, we have 
\begin{align}
\overline{\varphi}^{\mathcal{U}}_{\omega}(f) &= \varphi^{\mathcal{U}}_{\omega}(U_{\theta}f) \notag \\
&= \mathcal{U}\mathchar`-\lim_x \frac{1}{x} \int_0^x (U_{\theta}f)_{\omega}(t)dt \notag \\
&= \mathcal{U}\mathchar`-\lim_x \frac{1}{x} \int_0^x (U_{\theta}f_{\omega})(t)dt \notag \\
&= \mathcal{U}\mathchar`-\lim_x \frac{1}{x} \int_0^x \left(\frac{1}{\theta} \int_t^{t+\theta}f_{\omega}(s)ds\right)dt \notag \\
&= \mathcal{U}\mathchar`-\lim_x \frac{1}{x} \int_0^{x} \left(\frac{F_{\omega}(t+\theta) - F_{\omega}(t)}{\theta}) \right)dt, \notag 
\end{align}
where $F_{\omega}(x) = \int_0^x f_{\omega}(t)dt$. For every $x > 0$, we have by the dominated convergence theorem,
\[\lim_{\theta \to 0^+} \frac{1}{x} \int_0^x \left(\frac{F_{\omega}(t+\theta) - F_{\omega}(t)}{\theta}) \right)dt = \frac{1}{x} \int_0^x \lim_{\theta \to 0^+} \left(\frac{F_{\omega}(t+\theta) - F_{\omega}(t)}{\theta}) \right)dt = \frac{1}{x} \int_0^x f_{\omega}(t)dt. \]

Therefore, we get
\[\overline{\varphi}^{\mathcal{U}}_{\omega}(f) = \mathcal{U}\mathchar`-\lim_x \frac{1}{x} \int_0^x f_{\omega}(t)dt. \]
\end{prf}

We denote the set of all extensions of $\varphi_{\omega}^{\mathcal{U}} \in \mathcal{Q}$ by $\mathcal{Q}_1$. Then the following is immediate by Theorem 2.2 and Theorem 3.6.
\begin{thm}
$\mathcal{M}_1 = \overline{co}(\mathcal{Q}_1)$ holds.
\end{thm}

\section{Main results for multiplicative group}
Let $\mathbb{R}_+^{\times} = [1, \infty)$ and $L^{\infty}(\mathbb{R}_+^{\times})$ be the set of all essentially bounded measurable functions on $\mathbb{R}_+^{\times}$. Now for each $r \ge 1$ we introduce the following linear operator:
\[P_r : L^{\infty}(\mathbb{R}_+^{\times}) \longrightarrow L^{\infty}(\mathbb{R}_+^{\times}), \quad (P_rf)(x) = f(rx). \]
Let $P_r^*$ be its adjoint operator. Then we say that $\psi$ is a P-invariant mean if $\psi$ is a mean on $L^{\infty}(\mathbb{R}_+^{\times})$ and satifies
\[P_r^*\psi = \psi \quad for \ every \ r \ge 1. \]

Let us define a sublinear functional $\overline{L}_1$ on $L^{\infty}(\mathbb{R}_+^{\times})$ as
\[\overline{L}_1(f) = \lim_{\theta \to \infty} \limsup_{x \to \infty} \frac{1}{\log \theta} \int_{x}^{\theta x} f(t) \frac{dt}{t}. \]
We denote by $\mathcal{L}_1$ the class of means $\psi$ for which 
\[\psi(f) \le \overline{L}_1(f) \]
holds for every $f \in L^{\infty}(\mathbb{R}_+^{\times})$. 

Also let $\mathcal{M}$ be the class of means on $L^{\infty}(\mathbb{R}_+^{\times})$ for which
\[\varphi(f) \le \overline{M}(f) = \limsup_{x \to \infty} \frac{1}{x} \int_1^x f(t)dt \]
holds for every $f \in L^{\infty}(\mathbb{R}_+^{\times})$. Let us define the linear operator $U$ by
\[U : L^{\infty}(\mathbb{R}_+^{\times}) \longrightarrow L^{\infty}(\mathbb{R}_+^{\times}), \quad (Uf)(x) = \frac{1}{x} \int_1^x f(t)dt. \]
In particular, $\tilde{\mathcal{M}}$ be the subset of $\mathcal{M}$ consisting of those members $\varphi_{\omega}$ defined 
as follows.
\[\varphi_{\omega}(f) = e^{\omega}\mathchar`-\lim_x \frac{1}{x} \int_1^x f(t)dt, \]
where $\omega \in \Omega^*$ and $e^{\omega} = \{e^A : A \in \omega \}$ is an ultrafilter on $\mathbb{R}_+^{\times}$.

We define an algebraic isomorphism $W$ from $L^{\infty}(\mathbb{R}_+^{\times})$ onto $L^{\infty}(\mathbb{R}_+)$ as follows:
\[W : L^{\infty}(\mathbb{R}_+^{\times}) \longrightarrow L^{\infty}(\mathbb{R}_+), \quad (Wf)(x) = f(e^x). \]
Then we have the following commutative diagram
$$
\begin{CD}
L^{\infty}(\mathbb{R}_+^{\times}) @>W >> L^{\infty}(\mathbb{R}_+)\\
@VP_rVV @VVT_sV \\
L^{\infty}(\mathbb{R}_+^{\times}) @>W >> L^{\infty}(\mathbb{R}_+)
\end{CD}
$$
where $r=e^s, s \ge 0$. 

The relationship between the linear operators $S$ and $U$ can be given via $W$ as follows.

\begin{lem}
$U = W^{-1}SW$ holds.
\end{lem}

\begin{prf}
For any $f \in L^{\infty}(\mathbb{R}_+^{\times})$ and $x \ge 0$, we have
\[(SWf)(x) = \frac{1}{e^x} \int_0^x f(e^t)e^tdt = \frac{1}{e^x} \int_1^{e^x} f(s)ds.  \]
Hence
\[(W^{-1}SWf)(x) = \frac{1}{x} \int_1^x f(t)dt = (Uf)(x). \]
\end{prf}

Then it is easy to prove the following result.

\begin{thm}
$W^*\chi_{\omega} = \varphi_{\omega}$ holds for every $\omega \in \Omega^*$.
\end{thm}

\begin{prf}
For each $f \in L^{\infty}(\mathbb{R}_+)$, it holds that by Lemma 4.1,
\[(W^*\chi_{\omega})(f) = \omega\mathchar`-\lim_x (SWf)(x) = e^{\omega}\mathchar`-\lim_x (W^{-1}SWf)(x) = e^{\omega}\mathchar`-\lim_x (Uf)(x) = \varphi_{\omega}(f). \]
\end{prf}

\begin{thm}
$\mathcal{M}_1$ and $\mathcal{R}$ are affinely homeomorphic to $\mathcal{L}_1$ and $\mathcal{M}$ respectively via $W^*$, where $W^*$ is the adjoint operator of $W$.
\end{thm}

\begin{prf} 
We will show only for $\mathcal{M}_1$ and $\mathcal{L}_1$. The case of $\mathcal{R}$ and $\mathcal{M}$ can be proved similarly.
It is sufficient to prove that $\overline{L}_1(f) = \overline{M}_1(Wf)$ for every $f \in L^{\infty}(\mathbb{R}_+^{\times})$. 
\[
\limsup_{x \to \infty}\frac{1}{\theta} \int_x^{x + \theta} (Wf)(t)dt 
= \limsup_{x \to \infty} \frac{1}{\theta} \int_{e^x}^{e^x \cdot e^{\theta}} f(s) \frac{ds}{s}
= \limsup_{r \to \infty} \frac{1}{\theta} \int_r^{r \cdot e^{\theta}} f(s) \frac{ds}{s},   
\]
where we put $r=e^x$. And then put $y=e^{\theta}$ and we have 
\[\overline{M}_1(Wf) = \lim_{y \to \infty} \limsup_{r \to \infty} \frac{1}{\log y} \int_r^{ry} f(s) \frac{ds}{s} = \overline{L}_1(f). \]
\end{prf}

\medskip

Hence elements of $\mathcal{L}_1$ are P-invariant means since by the above diagram and elements of $\mathcal{M}_1$ are T-invariant means, for any $\psi \in \mathcal{L}_1$, let $\varphi = W^{*-1}\psi \in \mathcal{M}_1$ and then we have $\psi(P_rf) = (W^*\varphi)(P_rf) = \varphi(WP_rf) = \varphi(T_sWf) = \varphi(Wf) = (W^*\varphi)(f) = \psi(f)$.

\medskip
Concerning the class $\mathcal{M}$, in [5] we have shown the following results.

\begin{thm} 
$ex(\mathcal{C}) = \tilde{\mathcal{C}}$, where $ex(\mathcal{C})$ denotes the set of extreme points of $\mathcal{C}$.
\end{thm}

\begin{thm}
For any $\varphi \in \mathcal{C}$ there exists a unique probability measure $\mu$ on $\Omega^*$ such that
\[\varphi(f) = \int_{\Omega^*} \varphi_{\omega}(f) d\mu(\omega) \]
holds for every $f \in L^{\infty}(\mathbb{R}_+^{\times})$.
\end{thm}

Now it is obvious that Theorem 3.1 and Theorem 3.2 follows immediately from Theorem 4.3 and Theorem 4.4 with the aid of Theorem 4.1 and Theorem 4.2.

We define two subspaces $\Psi$ and $\Psi^{\prime}$ of $L^{\infty}(\mathbb{R}_+^{\times})$ by $\Psi = \{Uf : f \in L^{\infty}(\mathbb{R}_+^{\times})\}$ and $\Psi^{\prime} = \{f-Uf : L^{\infty}(\mathbb{R}_+^{\times})\}$, which are the counterparts of $\Phi$ and $\Phi^{\prime}$ respcetively defined in the former section. Then the following is obvious from Lemma 4.1.

\begin{lem}
$\Psi = W^{-1}\Phi$ and $\Psi^{\prime} = W^{-1}\Phi^{\prime}$ holds.
\end{lem}

Therefore, by Theorem 4.2 and Lemma 4.2, we get the following theorems immediately which correspond to Theorem 3.4 and Theorem 3.5 respectively.

\begin{thm}
For any mean $\psi$ on $L^{\infty}(\mathbb{R}_+^{\times})$, $\psi \in \mathcal{L}_1$ if and only if $\psi = 0$ on $\Psi^{\prime}$. Namely, $\psi \in \mathcal{L}_1$ if and only if $\psi$ is $U$-invariant.
\end{thm}

\begin{thm}
$\overline{L}_1(f) \le \overline{M}(f)$ for every $f \in L^{\infty}(\mathbb{R}_+^{\times})$.
\end{thm}

Given $f(x) \in L^{\infty}(\mathbb{R}_+^{\times})$, we consider the set of functions $\{f_r^{\times}(x)\}_{r \ge 1} \subseteq L^{\infty}(\mathbb{R}_+^{\times})$, where $f_r^{\times}(x) = f(rx)$. Then notice that this is a bounded set of $L^{\infty}(\mathbb{R}_+^{\times})$ and hence is a weak* relatively compact subset of $L^{\infty}(\mathbb{R}_+^{\times})$. Thus for any $\omega \in \Omega^*$ we can define its limit 
along $e^{\omega}$:
\[f_{\omega}^{\times}(x) = e^{\omega}\mathchar`-\lim_r f_r^{\times}(x). \]

Then we can also extend it to a function in $L^{\infty}(\mathbb{R}^{\times})$ in a similar way as $f_{\omega}$.

Now we take up the relation between $f_{\omega}(x)$ and $f_{\omega}^{\times}(x)$. For the sake of convenience, we define linear operators $T_{\omega}$ and $P_{\omega}$ by
\[T_{\omega} : L^{\infty}(\mathbb{R}_+) \longrightarrow L^{\infty}(\mathbb{R}), \quad (T_{\omega}f)(x) = f_{\omega}(x), \]
and
\[P_{\omega} : L^{\infty}(\mathbb{R}_+^{\times}) \longrightarrow L^{\infty}(\mathbb{R}^{\times}), \quad (P_{\omega}f)(x) = f_{\omega}^{\times}(x), \]
respectively. Then we have the following result (see [6] for a proof).

\begin{thm}
$P_{\omega} = W^{-1}T_{\omega}W$ for every $\omega \in \Omega^*$.
\end{thm}

Let $\omega \in \Omega^*$ and $\mathcal{U}$ be an ultrafilter on $\mathbb{R}_+^{\times}$ not containing any bounded set of $\mathbb{R}_+^{\times}$. Then we define for each $f \in L^{\infty}(\mathbb{R}_+^{\times})$
\[\psi^{\mathcal{U}}_{\omega}(f) = e^{\mathcal{U}}\mathchar`-\lim_x \frac{1}{\log x} \int_1^x \hat{f}_{\omega}(t)\frac{dt}{t}. \]
We denote by $\mathcal{P}$ the set of all such invariant means in $\mathcal{L}_1$.

\begin{thm}
$W^*\varphi_{\omega}^{\mathcal{U}} = \psi_{\omega}^{\mathcal{U}}$ holds.
\end{thm}

\begin{prf}
By Theorem 4.7 and integration by substitution, for every $f \in L^{\infty}(\mathbb{R}_+^{\times})$ we have
\begin{align}
(W^*\varphi^{\mathcal{U}}_{\omega})(f) = \varphi^{\mathcal{U}}_{\omega}(Wf) &= \mathcal{U}\mathchar`-\lim_x \frac{1}{x} \int_0^x (Wf)_{\omega}(t)dt \notag \\
&= \mathcal{U}\mathchar`-\lim_{x} \frac{1}{x} \int_0^x (Wf_{\omega}^{\times})(t)dt \notag \\
&= \mathcal{U}\mathchar`-\lim_x \frac{1}{x} \int_0^x f_{\omega}^{\times}(e^t)dt \notag \\
&= \mathcal{U}\mathchar`-\lim_x \frac{1}{x} \int_1^{e^x} f_{\omega}^{\times}(t)\frac{dt}{t} \notag \\
&= e^{\mathcal{U}}\mathchar`-\lim_x \frac{1}{\log x} \int_1^x f_{\omega}^{\times}(t)\frac{dt}{t} \notag \\
&= \psi^{\mathcal{U}}_{\omega}(f). \notag 
\end{align}

\end{prf}

In particular, by Theorem 3.8, Theorem 4.2 and Theorem 4.8, we obtain the following result analogous to Theorem 3.8.
\begin{thm}
$\mathcal{L}_1 = \overline{co}(\mathcal{P})$.
\end{thm}

\section{Applications to summability methods}
Recall that for a function $f(n)$ on $\mathbb{N}$, its Ces\`{a}ro mean $M_d(f)$ is defined as 
\[M_d(f) = \lim_{n \to \infty} \frac{1}{n} \sum_{i=1}^n f(i) \]
if the limit exists. It can be seen as a summability method of the most simple type. Naturally, for a measurable function $f(x)$ on $\mathbb{R}_+^{\times}$, an integral version of Ces\`{a}ro mean is defined as
\[M(f) = \lim_{x \to \infty} \frac{1}{x} \int_1^x f(t)dt \]
if the limit exists. Notice that restricting the domain of $M$ to bounded measurable functions on $\mathbb{R}_+^{\times}$,   $M$ can be viewed as a continuous linear functional on the subspace $\mathcal{D}(M)$ of $L^{\infty}(\mathbb{R}_+^{\times})$ whose elements possess the above limit.

In this section we will study summability methods by functional analytic methods. In what follows, we will restrict ourselves only to bounded measurable functions on $\mathbb{R}_+^{\times}$ and summability methods mean pairs $(F, \mathcal{D}(F))$ of the domain $\mathcal{D}(F)$ of $F$, which is a closed subspace of $L^{\infty}(\mathbb{R}_+^{\times})$, and a continuous linear functional $F : \mathcal{D}(F) \rightarrow \mathbb{R}$.

In particular, the summability method $M$ can be formulated in terms of the class of linear functionals $\mathcal{M}$: we define for each $f \in L^{\infty}(\mathbb{R}_+^{\times})$
\[\overline{M}(f) = \sup_{\varphi \in \mathcal{M}} \varphi(f) = \limsup_{x \to \infty} \frac{1}{x} \int_1^x f(t)dt \]
and we also define the lower version of $\overline{M}$ by
\[\underline{M}(f) := -\overline{M}(-f) = \inf_{\varphi \in \mathcal{M}} \varphi(f) = \liminf_{x \to \infty} \frac{1}{x} \int_1^x f(t)dt. \]
Then since we have 
\[M(f) = \alpha \Longleftrightarrow \overline{M}(f) = \underline{M}(f) = \alpha, \]
it holds that 
\[M(f) = \alpha \Longleftrightarrow \varphi(f) = \alpha \ for \ every \ \varphi \in \mathcal{M}. \]

Now we consider one of the generalizeations of $M$, which is an integral version of summability methods introduced by H\"{o}lder as generalizations of Ces\`{a}ro mean, namely, iterations of $M$ (see [3] for details). Notice that for $f \in L^{\infty}(\mathbb{R}_+^{\times})$ it can be written as
$M(f) = \lim_{x \to \infty} (Uf)(x)$. Then we define $H_2(f)$ by
\[H_2(f) = \lim_{x \to \infty} \frac{1}{x} \int_1^x (Uf)(t)dt = \lim_{x \to \infty} \frac{1}{x} \int_1^x \left(\frac{1}{t} \int_1^t f(t)dt\right)dx = \lim_{x \to \infty} (U^2f)(x), \]
if this limit exists.
We can repeat this procedure inductively and get the sequence of summability methods $(M, \mathcal{D}(M)) = (H_1, \mathcal{D}(H_1)), (H_2, \mathcal{D}(H_2)), \ldots, (H_k, \mathcal{D}(H_k)), \ldots$, where $H_k$ is defined by
\[H_k(f) = H_{k-1}(Uf) = \lim_{x \to \infty} (U^kf)(x), \quad k=1, 2, \ldots . \]

Similarly, we can relate each $(H_k, \mathcal{D}(H_k))$ to a sublinear functional $\overline{H}_k$ defined as 
\[\overline{H}_k(f) = \limsup_{x \to \infty} (U^kf)(x), \quad f(x) \in L^{\infty}(\mathbb{R}_+^{\times}), \]
or to a weak* compact convex subset $\mathcal{H}_k$ of $L^{\infty}(\mathbb{R}_+^{\times})^*$ whose elements $\varphi$ satisfy the condition that $\varphi(f) \le \overline{H}_k(f)$ holds for every $f \in L^{\infty}(\mathbb{R}_+^{\times})$.

Moreover, notice that the sequence $\{\overline{H}_k\}_{k=1}^{\infty}$ of sublinear functionals is monotonically decreasing:
\[\overline{H}_1(f) \ge \overline{H}_2(f) \ge \ldots \overline{H}_k(f) \ge \ldots \]
for every $f \in L^{\infty}(\mathbb{R}_+^{\times})$. This is obviously bounded below and there exists a limit $\overline{H}_{\infty}(f) := \lim_{k \to \infty} \overline{H}_k(f)$ for each $f \in L^{\infty}(\mathbb{R}_+^{\times})$. It is easy to see that this functional $\overline{H}_{\infty} : L^{\infty}(\mathbb{R}_+^{\times}) \rightarrow \mathbb{R}$ is also sublinear and it defines a summability method $(H_{\infty}, \mathcal{D}(H_{\infty}))$ by $H_{\infty}(f) = \alpha$ if and only if $\overline{H}_{\infty}(f) = \underline{H}_{\infty}(f) = \alpha$ or by $H_{\infty}(f) = \alpha$ if and only if $\varphi(f) = \alpha$ for every $\varphi \in \mathcal{H}_{\infty}$, where $\underline{H}_{\infty}(f) = -\overline{H}_{\infty}(-f)$ and $\mathcal{H_{\infty}}$ is a weak* compact convex subset of $L^{\infty}(\mathbb{R}_+^{\times})^*$ whose elements $\varphi$ satisfy the condition that $\varphi(f) \le \overline{H}_{\infty}(f)$ holds for every $f \in L^{\infty}(\mathbb{R}_+^{\times})$. Then it is noted that $\mathcal{H}_{\infty} = \cap_{k=1}^{\infty} \mathcal{H}_k$ holds. Our main aim of this section is to show the following theorem.

\begin{thm}
$\overline{H}_{\infty}(f) = \overline{L}_1(f)$ holds for every $f \in L^{\infty}(\mathbb{R}_+^{\times})$. In particular, $(H_{\infty}, \mathcal{D}(H_{\infty}))$ and $(L_1, \mathcal{D}(L_1))$ are the same summability method.
\end{thm}

Since it is convenient to work in the setting of Section 3, in the following, we will first formulate and prove the corresponding theorem in the context of Section 3 and then transfer it to the above theorem, as in the Section 4, via isomorphisms $W$ and $W^*$.

Now let $(R, \mathcal{D}(R))$ be the summability method defined on the subspace $\mathcal{D}(R)$ of $L^{\infty}(\mathbb{R}_+)$ whose elements $f$ have the limit
\[R(f) = \lim_{x \to \infty} \frac{1}{e^x} \int_0^x f(t)e^tdt. \]
In the same way as above, we can consider its iterations $E_1=R, E_2, E_3, \ldots$, defined by
\[E_k : \mathcal{D}(E_k) \longrightarrow \mathbb{R}, \quad E_k(f) = \lim_{x \to \infty} (S^kf)(x), \quad k =1, 2, \ldots. \]

Also, we define sublinear functionals $\overline{E}_k$ by
\[\overline{E}_k(f) = \limsup_{x \to \infty} (S^kf)(x), \quad k=1,2, \ldots. \]
Let us $\overline{E}_{\infty}(f) = \lim_{k \to \infty} \overline{E}_k(f)$ and $(E_{\infty}, \mathcal{D}(E_{\infty}))$ denotes the induced summability method defined as above. Also for $k = 1,2, \ldots, \infty$, we denote by $\mathcal{E}_k$ the set of linear functionals $\varphi$ on $L^{\infty}(\mathbb{R}_+)$ such that $\varphi(f) \le \overline{E}_k(f)$ holds for every $f \in L^{\infty}(\mathbb{R}_+)$. Then it also holds that $\mathcal{E}_{\infty} = \cap_{k=1}^{\infty} \mathcal{E}_k$.
Now we can formulate a version of Theorem 5.1 as follows.

\begin{thm}
$\overline{E}_{\infty}(f) = \overline{M}_1(f)$ holds for every $f \in L^{\infty}(\mathbb{R}_+)$. In particular, $(E_{\infty}, \mathcal{D}(R_{\infty}))$ and $(M_1, \mathcal{D}(M_1))$ are the same summability method.
\end{thm}

Recall that the convolution $(f*\phi)(x) \in C_{ub}(\mathbb{R})$ of $f(x) \in L^{\infty}(\mathbb{R})$ and $\phi(x) \in L^1(\mathbb{R})$ is defined by 
\[(f*\phi)(x) = \int_{\mathbb{R}} f(x-t)\phi(t)dt, \quad x \in \mathbb{R}. \]

\begin{lem}
For each $f \in L^{\infty}(\mathbb{R}_+)$ and $\omega \in \Omega^*$, $(Sf)_{\omega}(x) = \int_0^{\infty} f_{\omega}(x-t)e^{-t}dt = (f_{\omega} * h)(x)$ holds, where $h$ is a function in $L^1(\mathbb{R})$ such that
\[h(x) = \begin{cases}
            e^{-x} & \text{if $x \ge 0$}, \\
            0       & \text{if $x < 0$}.
         \end{cases}     
\]
\end{lem}

\begin{prf}
By the definition of $(Sf)_{\omega}(x)$ we have 
\[(Sf)_{\omega}(x) = \tau^x \omega\mathchar`-\lim_s \frac{1}{e^s} \int_0^s f(t)e^tdt = \int_0^{\infty} f_{\tau^x\omega}(-t)e^{-t}dt = \int_0^{\infty} f_{\omega}(x-t)e^{-t}dt = (f_{\omega} * h)(x). \]
\end{prf}

Let us denote $h^n = \overbrace{h * \ldots * h}^n$ for $n \ge 1$. Then we have the following lemma, which can be proved by a direct computation and we omit the proof.
\begin{lem}
For $n = 1,2, \ldots$, it holds that
\[h^n(x) = \begin{cases}
                e^{-x} \cdot \frac{x^{n-1}}{(n-1)!} & \text{if $x \ge 0$}, \\
                0 & \text{if $x < 0$}.
              \end{cases}
\]
\end{lem}

\begin{thm}
For each $f \in L^{\infty}(\mathbb{R}_+), \omega \in \Omega^*$ and $n \ge 1$, we have 
\[(S^nf)_{\omega}(x) = (f_{\omega} * h^n)(x) = \int_0^{\infty} f_{\omega}(x-t)e^{-t} \frac{t^{n-1}}{(n-1)!}dt. \]
\end{thm}

\begin{prf}
By Lemma 5.1 and Lemma 5.2, it is easy to see that 
\begin{align}
(S^2f)_{\omega}(x) &= (S(Sf))_{\omega}(x) = ((Sf)_{\omega} * h)(x) = ((f_{\omega} * h) * h)(x)  \notag \\
&= (f_{\omega} * (h * h))(x) = (f_{\omega} * h^2)(x) = \int_0^{\infty} f_{\omega}(x-t)e^{-t} tdt. \notag
\end{align}
For a general $n \ge 1$, we get the result by induction.
\end{prf}

\begin{lem}
For each $f \in L^{\infty}(\mathbb{R}_+), \omega \in \Omega^*$ and $n \ge 1$, we have 
\[|\overline{(S^nf)}(\omega) - \overline{(S^{n+1}f)}(\omega)| \le 2\frac{e^{-n} \cdot n^n}{n!} \dot \|f\|_{\infty}. \]
\end{lem}

\begin{prf}
It holds that
\begin{align}
|\overline{(S^nf)}(\omega) - \overline{(S^{n+1}f)}(\omega)| &=  
\left|\int_0^{\infty} f_{\omega}(-t)e^{-t} \frac{t^{n-1}}{(n-1)!}dt - \int_0^{\infty} f_{\omega}(-t)e^{-t} \frac{t^n}{n!}dt \right| \notag \\
&\le \int_0^{\infty} |f_{\omega}(-t)| e^{-t} \left|\frac{t^{n-1}}{(n-1)!} - \frac{t^n}{n!}\right|dt \notag \\
&\le \|f_{\omega}\|_{\infty} \int_0^{\infty} e^{-t} \left|\frac{t^{n-1}}{(n-1)!} - \frac{t^n}{n!}\right|dt \notag \\
&= \|f_{\omega}\|_{\infty} \int_0^n e^{-t} \left(\frac{t^{n-1}}{(n-1)!} - \frac{t^n}{n!}\right)dt \notag \\
&+ \|f_{\omega}\|_{\infty} \int_n^{\infty} e^{-t} \left(\frac{t^n}{n!} - \frac{t^{n-1}}{(n-1)!}\right)dt \notag \\
&= \|f_{\omega}\|_{\infty} \left(\left[\frac{e^{-t}t^n}{n!}\right]_0^n - \left[\frac{e^{-t}t^n}{n!}\right]_n^{\infty}\right) 
= 2\|f_{\omega}\|_{\infty} \frac{e^{-n} n^n}{n!}. \notag
\end{align}
\end{prf}

Now we can prove Theorem 5.2. 
\begin{prf}[Proof of Theorem 5.2]
Since, as we have seen in the proof of Theorem 3.4, $\overline{M}_1(f) = \overline{M}_1(Sf)$, we have that for all $n \ge 1$
\[\overline{M}_1(f) = \overline{M}_1(S^nf) \le \limsup_{x \to \infty} (S^nf)(x) = \overline{E}_n(f).  \]
Then letting $n \to \infty$, we get that
\[\overline{M}_1(f) \le \lim_{n \to \infty} \overline{E}_n(f) = \overline{E}_{\infty}(f). \]
On the other hand, by Lemma 5.3 together with Stirling's formula $n! \sim \sqrt{2\pi n}(\frac{n}{e})^n$ it holds that
\begin{align}
|\overline{E}_{\infty}(f-Sf)| &= \lim_{n \to \infty} |\overline{E}_n(f-Sf)| = \lim_{n \to \infty} \limsup_{x \to \infty} |(S^nf)(x) - (S^{n+1}f)(x)| \notag \\ 
&= \lim_{n \to \infty} \sup_{\omega \in \Omega^*} |\overline{(S^nf)}(\omega) - \overline{(S^{n+1}f)}(\omega)| \notag \\
&\le \lim_{n \to \infty} 2\frac{e^{-n} \cdot n^n}{n!} \dot \|f\|_{\infty} = 0. \notag
\end{align}
Therefore, we have that if $\varphi \in \mathcal{E}_{\infty}$ then $\varphi = 0$ on $\Phi^{\prime}$, which means that $\varphi \in \mathcal{M}_1$ by Theorem 3.4. Thus 
\[\overline{E}_{\infty}(f) = \sup_{\varphi \in \mathcal{E}_{\infty}} \varphi(f) \le \sup_{\varphi \in \mathcal{M}_1} \varphi(f) = \overline{M}_1(f). \]
Hence we obtain that
\[\overline{M}_1(f) = \overline{E}_{\infty}(f). \]
\end{prf}

\begin{prf}[Proof of Theorem 5.1]
First of all, notice that
\[W^* : \mathcal{E}_i \longrightarrow \mathcal{H}_i \]
is an affine homeomorphism between $\mathcal{E}_i$ and $\mathcal{H}_i$ for every $i \ge 1$. In fact, by Lemma 4.1, for each $i \ge 1$ we have 
\begin{align}
\overline{H}_i(f) &= \limsup_{x \to \infty}  (U^if)(x) \notag \\
&= \limsup_{x \to \infty} (W^{-1}S^iWf)(x) \notag \\
&= \limsup_{x \to \infty} (S^iWf)(\log x) \notag \\
&= \limsup_{x \to \infty} (S^iWf)(x) \notag \\
&= \overline{E}_i(Wf). \notag
\end{align}
Recall that $W^* \mathcal{M}_1 = \mathcal{L}_1$ by Theorem 4.2. Also by Theorem 5.2, $\mathcal{M}_1 = \cap_{i=1}^{\infty} \mathcal{E}_i$ holds and these implies that $\mathcal{L}_1 = \cap_{i=1}^{\infty} \mathcal{H}_i$, which means that $\overline{L}_1(f) = \overline{M}_{\infty}(f)$. This completes the proof.
\end{prf}

\end{document}